\documentclass[12pt,a4paper]{article}
\usepackage{amsmath,amssymb,amsthm,xspace}
\usepackage{amscd,latexsym,amscd,stmaryrd,xy}
\xyoption{all} 
\usepackage{enumerate}

\newtheorem{lemma}{Lemma}[section]
\newtheorem{remark}{Remark}[section]
\newtheorem{theorem}{Theorem}[section]

\newtheorem{proposition}{Proposition}[section]

\newtheorem{corollary}{Corollary}[section]

\newcommand{\C}{\ensuremath{\mathbb{C}}}

\newcommand{\Z}{\ensuremath{\mathbb{Z}}}

\font\sc=rsfs10 at 12 pt
\font\scs=rsfs10 at 10 pt
\font\scb=rsfs10 at 16 pt
\font\scbb=rsfs10 at 18 pt

\newcommand{\cT}{\sc\mbox{T}\hspace{1.0pt}}
\newcommand{\cP}{\sc\mbox{P}\hspace{1.0pt}}
\newcommand{\OO}{\sc\mbox{O}\hspace{1.0pt}(\mathfrak{p},\Lambda)}
\newcommand{\cO}{\sc\mbox{O}\hspace{1.0pt}}
\newcommand{\ccO}{\scb\mbox{O}\hspace{1.0pt}}

\newcommand{\ccccO}{\scbb\mbox{O}\hspace{1.0pt}}
\newcommand{\cX}{\sc\mbox{X}\hspace{1.0pt}}
\newcommand{\ccX}{\scs\mbox{X}}

\newcommand{\End}{\operatorname{End}}
\newcommand{\Hom}{\operatorname{Hom}}
\renewcommand{\hom}{\operatorname{hom}}
\newcommand{\Ext}{\operatorname{Ext}}
\newcommand{\ext}{\operatorname{ext}}
\newcommand{\tto}{\twoheadrightarrow}

\newcommand{\fg}{\mathfrak{g}}
\newcommand{\fn}{\mathfrak{n}}
\newcommand{\fh}{\mathfrak{h}}
\newcommand{\bfW}{\mathbf{W}}
\newcommand{\bfV}{\mathbf{V}}
\newcommand{\bfG}{\mathbf{G}}
\newcommand{\bfH}{\mathbf{H}}

\newcommand{\bfC}{\mathbf{C}}
\newcommand{\ttA}{\mathtt{A}}
\newcommand{\ttB}{\mathtt{B}}
\newcommand{\ttC}{\mathtt{C}}
\newcommand{\ttD}{\mathtt{D}}

\begin{document}
\title{Applications of the category of linear complexes of tilting modules 
associated with the category $\ccccO$}
\author{Volodymyr Mazorchuk}
\date{}

\maketitle
\begin{abstract}
We use the category of linear complexes of tilting modules for the BGG 
category $\cO$, associated with a semi-simple complex finite-dimensional
Lie algebra $\fg$, to reprove in purely algebraic way several known 
results about $\cO$ obtained earlier by different authors using geometric 
methods. We also obtain several new results about the parabolic category 
$\OO$.
\end{abstract}

\section{Introduction and preliminaries}\label{s1}

Let $\fg$ be a semi-simple complex finite-dimensional Lie algebra with a
fixed triangular decomposition, $\fg=\fn_-\oplus\fh\oplus\fn_+$, and
$\pi$ be the corresponding basis of the root system of $\fg$. 
Let further $\rho$ be the half of the sum of all positive roots and
$\bfW$ be the Weyl group of $\fg$. For $w\in\bfW$ and $\lambda\in\fh^*$ set 
$w\cdot\lambda=w(\lambda+\rho)-\rho$. Denote by $w_0$ the longest element of 
$\bfW$ and by $\mathfrak{l}$ the length function on $\bfW$. We consider
$\bfW$ as the partially ordered set with respect to the Bruhat order $<$ under 
the convention that the identity $e$ is the smallest element.

Consider the principal block  $\cO_0$ of the BGG category $\cO$ for $\fg$, 
associated with the  triangular decomposition above (see \cite{BGG1,So1}). 
This category is a highest weight category in the sense of \cite{CPS} and hence 
is equivalent to the module category of some finite-dimensional associative
quasi-hereditary basic algebra, which we denote by $\ttA$ (see \cite{BGG1,DR}).
The Chevalley involution on $\fg$ induces a natural involutive contravariant
exact self-equivalence on $\cO$, and hence on $\ttA\mathrm{-mod}$, 
which we will denote by $\star$. The isomorphism classes of simple modules in 
$\cO_0$ (and thus also in $\ttA\mathrm{-mod}$) are in natural bijection with 
the elements of $\bfW$ under the convention that $e$ corresponds to the 
one-dimensional module in $\cO_0$. For $w\in \bfW$ we introduce the following 
notation:
\begin{itemize}
\item $L(w)$ is the simple $\ttA$-module, which corresponds to the simple
highest weight module in $\cO$ with the highest weight $w\cdot 0$ 
(see \cite[Chapter~7]{Di});
\item $\Delta(w)$ is the standard $\ttA$-module corresponding to $L(w)$ (in 
$\cO$ the module $\Delta(w)$ is the Verma module with the highest weight 
$w\cdot 0$, see \cite[Chapter~7]{Di});
\item $\nabla(w)=\Delta(w)^{\star}$ is the corresponding costandard module
for $\ttA$ (the dual Verma module in $\cO$);
\item $P(w)$ is the projective cover of $L(w)$;
\item $I(w)=P(w)^{\star}$ is the injective envelope of $L(w)$;
\item $T(w)$ is the indecomposable tilting $\ttA$-module, which corresponds
to $\Delta(w)$ (see \cite{Ri1}).
\end{itemize}
The principal result about the category $\cO$ is the so-called {\em Kazhdan-Lusztig} 
(or, simply, the KL-) Theorem, proved in \cite{BeBe,BrKa}, which describes the 
composition multiplicities of $\Delta(w)$ in terms of the Kazhdan-Lusztig combinatorics 
of the Hecke algebra, see \cite{KaLu}. This result and the geometric approach to 
$\cO$ have been  used by Soergel in \cite{So1} to show that $\ttA$ is Koszul and even 
Koszul self-dual. This was further extended in \cite{BeGiSo}, where it was shown that 
all associative algebras associated with the blocks of $\cO$ are Koszul and that the 
Koszul dual of the algebra of a singular block of $\cO$ is the algebra of the regular 
block of certain parabolic generalization of $\cO$, introduced by Rocha-Caridi in 
\cite{RoCa}. In \cite{Ba} this result  was  even further extended to all 
parabolic-singular blocks. Both \cite{BeGiSo} and \cite{Ba} use ``heavy'' geometric 
arguments.

The algebra $\ttA$ can be given a nice combinatorial description via the coinvariant
algebra $\bfC$ associated with $\bfW$. In \cite{So1} it is shown that 
$\bfC\cong\End_{\ttA}\left(P(w_0)\right)$; that the module $P(w_0)$, viewed as
a $\bfC$-module, can be described combinatorially; and that $\ttA\cong
\End_{\bfC}\left(P(w_0)_{\bfC}\right)$ (see also \cite{KoSlXi} for the last isomorphism). 
Furthermore, in \cite{So1} it is shown that $\End_{\ttA}\left(P(w_0)\right)$ is 
isomorphic to the quotient of the polynomial algebra over a homogeneous ideal, in 
particular, that it is $\Z$-graded (however the grading is not unique). If we let the 
generators of $\End_{\ttA}\left(P(w_0)\right)$ to have degree $2$, we obtain the grading  
which coincides with the grading on $\bfC$ obtained via interpretation of $\bfC$ as 
the cohomology algebra of a certain flag manifold (see for example \cite{Hi}). Further, 
the combinatorial construction of Soergel implies that $P(w_0)$ is a $\Z$-graded 
$\bfC$-module, which makes $\ttA$ into a $\Z$-graded algebra. In \cite{BeGiSo} it was 
shown that this grading is the Koszul grading on $\ttA$ using some arguments from the 
``mixed'' geometry.

The paper \cite{St1} initiates the study of the graded version of the category $\cO$
with respect to the ``natural'' grading described above. It is important that
Soergel's combinatorial description and hence the natural grading on $\ttA$ do not 
depend on the KL-Theorem. In particular, in \cite{St1} the author obtains
several results about the graded version of the category $\cO$, which do not
depend on the KL-Theorem. It seems that there is a hope that the graded approach
might be a way to give an algebraic proof of the KL-Theorem. Because of Vogan's
formulation of the KL-Theorem, see \cite{Vo}, and the graded description of the translation
functors from \cite{St1}, the KL-Theorem would follow if one would show that the natural
grading on $\ttA$ is positive in the sense that all non-zero components have
non-negative degrees, and the radical of $\ttA$ coincides with the sum of all
components of positive degrees. However, so far this seems to be very difficult.

Let us now consider $\ttA$ as a graded algebra with respect to the natural grading. 
Then it is easy to see that all simple, projective, injective, standard and 
costandard modules admit graded lifts. In \cite{MaOv,Zh} it was shown that in this 
case all tilting modules admit graded lifts as well. In particular, the Ringel dual 
$R(\ttA)$ of $\ttA$ (see \cite{Ri1}) is automatically graded. Fixing natural graded lifts of
indecomposable tilting modules one can consider the category $\cT(\ttA)$ of linear 
complexes of tilting modules for $\ttA$. In \cite{MaOv} it was shown that $\cT(\ttA)$ 
is equivalent to the category of locally finite-dimensional  graded modules over the 
quadratic dual $R(\ttA)^{!}$ of $R(\ttA)$. Using the Koszul self-duality of $\ttA$ (see 
\cite{So1}) and the Ringel self-duality of $\ttA$ (see \cite{So2}), one obtains that 
$\cT(\ttA)$ is equivalent to the category of locally finite-dimensional graded 
$\ttA$-modules. In the present paper we would like to derive several consequences from 
this fact.

In the present paper we use two ``heavy'' prerequisites. The first one is the KL-Theorem.
Using the tilting module approach together with the KL-combinatorics (which follows from 
the KL-Theorem) we show in Section~\ref{s2} that $\ttA$ is Koszul and that the natural 
grading on $\ttA$ is the Koszul one. However, we are not able to derive the Koszul
self-duality for $\ttA$ by our methods, hence we use this result of Soergel as the
second ``heavy'' prerequisite. In Section~\ref{s3} we give a brief synopsis about the 
category of linear complexes of tilting modules for a quasi-hereditary algebra, studied in 
\cite{MaOv}, and list some corollaries for the category $\cO_0$. In Section~\ref{s4} we 
show that the associative algebras, associated to all singular blocks of $\cO$ are Koszul. 
Already on this stage we need Koszul self-duality for $\ttA$. We also use the machinery
developed in \cite{MaOv}. In Section~\ref{s5} we recall the definition of the
parabolic category $\cO_S$ due to Rocha-Caridi, \cite{RoCa}, and show that all
(regular and singular) blocks of this category are Koszul. This result happens to be
a relatively easy corollary from the corresponding result for the singular blocks of 
$\cO$, which we prove earlier. It is possibly interesting to point out that the statement 
that the regular block of $\cO_S$ is Koszul does not require Soergel's Koszul
self-duality result. In Section~\ref{s6} we define a new properly stratified (in the sense 
of \cite{Dl}) parabolic generalization of $\cO$ and $\cO_S$. Roughly speaking, it is 
the translation of a singular block of $\cO$ out of the wall. We use this category 
and the corresponding category of linear complexes of tilting modules to reprove in 
Section~\ref{s7} the parabolic-singular Koszul dualities from \cite{BeGiSo,Ba} in a 
purely algebraic way. 

Apart from the KL-Theorem and Soergel's Koszul self-duality result, the main 
ingredients of this proof are:
\begin{itemize}
\item the way the composition of the Ringel and Koszul self-dualities on $\ttA$ 
permute the primitive idempotents, see  \cite{So1,So2};
\item the way the quadratic duality, defined via linear complexes of tilting modules, 
works.
\end{itemize}

The simple modules in a singular block of $\cO$ are indexed by certain left cosets in 
$\bfW$. The simple modules in the corresponding parabolic category are indexed by the 
corresponding right cosets. The composition of the Ringel and Koszul self-dualities of 
$\ttA$ switches the left and the right cosets in $\bfW$ and, applied to a certain 
subcategory of $\cT(\ttA)$, gives the necessary duality. It is again perhaps 
interesting to note that for a finite-dimensional semi-simple Lie algebra $\fg$ 
the Ringel self-duality of $\cO_0$ does not require the approach proposed in 
\cite{So2}, which is based on Arkhipov's twisting functor. Instead one can use a 
direct approach from \cite{FuKoMa}, based on the translation functors, which 
substantially simplifies the argument. 

The necessary preliminaries about the quadratic duality are collected in Section~\ref{s3},
where it is shown that this duality always switches the centralizer algebras, similar to 
those considered in Section~\ref{s6}, and the quotient algebras, similar to the blocks of 
$\cO_S$.

Finally, we use $\cT(\ttA)$ to study the properties of the parabolic generalization
$\OO$ of $\cO$, considered in \cite{FuKoMa} (or, equivalently, using \cite{BeGe},
certain singular blocks in the category of Harish-Chandra bimodules for $\fg$).
The blocks of this category are equivalent to the module categories of certain properly 
stratified algebras, see \cite{FuKoMa}. Most of our results here are about the
regular blocks. We show that the standard modules in the regular blocks 
have linear tilting coresolutions and that the costandard  modules in the regular 
blocks have linear tilting resolutions. This can be proved in at least two different ways.
The first way is analogous to that for $\cO_0$ and uses the properties of translation 
functors. The second way uses the existence of a non-standard BGG-type resolutions for 
certain highest weight modules. We further compute the quadratic dual to the properly 
stratified algebra of the regular block of $\OO$, which happens to be the parabolic 
quotient of $\ttA$ similar to $\cO_S$, but related to the left cosets instead of the 
right cosets. This is one more evidence for the strong asymmetry of $\cO$ with respect 
to the left and the right cosets of $\bfW$ (the first one was obtained in \cite{MaSt}). 
However, this result gives a possibility to describe the lawyers of the tilting 
(co)resolutions of standard and costandard modules from $\OO$ in terms of the 
Kazhdan-Lusztig combinatorics. We also derive several facts about the extensions 
between standard and proper standard modules in $\OO$.

\section{Koszulity of the natural grading for the regular blocks}\label{s2}

In this section we use only the KL-Theorem (in particular, in Vogan's formulation)
and do not use Soergel's Koszul self-duality of $\ttA$. We recall that a 
quasi-hereditary algebra is called {\em standard Koszul}, \cite{ADL}, provided 
that it is positively graded, all standard modules admit linear projective resolutions 
(meaning the the $l$-th term of the resolution is generated in degree $l$), and all 
costandard modules admit linear injective coresolutions (meaning the the
$l$-th term of the coresolution is cogenerated in degree $-l$). By \cite[Theorem~1]{ADL},
every standard Koszul quasi-hereditary algebra is Koszul.

In all categories of graded modules and for all $k\in\Z$ we denote by $\langle k\rangle$ 
the functor of the shift of grading, which maps the degree $l$ to the degree $l-k$.
For categories of complexes and for all $k\in\Z$ we denote by $[k]$ the functor
of the shift of the position in a complex, which maps the position $l$ to 
the position $l-k$.

According to \cite{So1} and \cite{KoSlXi}, the $\bfC$-module $P(w_0)$ admits a 
decomposition, $P(w_0)_{\bfC}=\oplus_{w\in\bfW}D(w)$, with indecomposable modules $D(w)$
constructed recursively as follows. $D(e)$ is the simple $\bfC$-module,
which we consider as the graded module, concentrated in degree $0$.
For a simple reflection, $s\in\bfW$, we denote by $\bfC^s$ the subalgebra of
$s$-invariants in $\bfC$. Let $w\in\bfW$ and $w=s_1\dots s_k$ be a reduced 
decomposition of $w$. Then the module 
\begin{displaymath}
\tilde{D}(w)=\bfC\otimes_{\bfC^{s_k}}\bfC\otimes_{\bfC^{s_{k-1}}}\dots
\otimes_{\bfC^{s_1}}D(e)\langle k\rangle
\end{displaymath}
has one dimensional component of degree $-k$. The module $D(w)$ is the
indecomposable direct summand of $\tilde{D}(w)$ such that $D(w)_{-k}\neq 0$.
This fixes a grading on $P(w_0)_{\bfC}$ and makes the algebra 
$\ttA=\End_{\bfC}\left(P(w_0)_{\bfC}\right)$ into a graded algebra, 
see \cite{So1}. We call this grading on $\ttA$ {\em natural}. We denote by
$\ttA\mathrm{-mod}$ and $\ttA\mathrm{-grmod}$ the categories of all finitely
generated $\ttA$-modules and all finitely generated graded (with respect to
the natural grading) $\ttA$-modules respectively. Remark that the morphisms in
$\ttA\mathrm{-grmod}$ are homogeneous maps of degree $0$. We set
$\ext_{\ttA}=\Ext_{\ttA\mathrm{-grmod}}$, $\hom_{\ttA}=\Hom_{\ttA\mathrm{-grmod}}$.

For standard graded lifts we will use the same symbol as for ungraded modules.
In particular, we set $L=\oplus_{w\in \bfW}L(w)$ and same for $P$, $I$, $T$, 
$\Delta$, $\nabla$, both for graded and ungraded modules.
We concentrate $L$ in degree $0$ and fix a grading on $P$ such that the
natural map $P\tto L$ is a morphism in $\ttA\mathrm{-grmod}$. Further, we 
fix a grading on $I$ such that the natural map $L\hookrightarrow I$ is a morphism 
in $\ttA\mathrm{-grmod}$. The natural maps $P\tto\Delta$ and $\nabla\tto I$ then
automatically induce gradings on $\Delta$ and $\nabla$. Further, we fix the
grading on $T$ such that the natural map $\Delta\hookrightarrow T$ is
a morphism in $\ttA\mathrm{-grmod}$. It follows then automatically that
the natural map $T\tto\nabla$ is a morphism in $\ttA\mathrm{-grmod}$, we
refer the reader to \cite{MaOv} for details. 

By \cite[Section~6]{St1} the duality $\star$ lifts to a duality on 
$\ttA\mathrm{-grmod}$, which we will denote by the same symbol. Note 
that $\star$ acts on degrees via multiplication with $-1$.
Let $s$ be a simple reflection and $\theta_s$ be the translation functor
through the $s$-wall, see \cite{Ja}. This functor is exact and self-adjoint.
In \cite{St1} it was shown that $\theta_s$ admits a graded lift, that is
it lifts to an exact and self-adjoint functor on $\ttA\mathrm{-grmod}$.
Furthermore, Vogan's version of the KL-Theorem, \cite{Vo}, asserts that this
lift can be chosen such that $(\theta_s L)_i=0$ for all $i\neq -1,0,1$ and
such that $(\theta_s L)_0$ is a semi-simple module.  Now we are going to use 
this to prove the following main result of the present section, first proved in 
\cite[Theorem~18]{So1} and \cite[4.5]{BeGiSo}.

\begin{theorem}\label{t1}
$\ttA$ is a standard Koszul quasi-hereditary algebra and the natural grading on 
$\ttA$ is Koszul.
\end{theorem}

To prove this we will need some preparation. We start with the following result,
which is a graded version of \cite[Proposition~4]{FuKoMa}.

\begin{proposition}\label{p2}
\begin{enumerate}[(i)]
\item\label{p2.1} Let $w\in\bfW$ then there is an inclusion (of graded modules)
$T(w_0w)\langle-\mathfrak{l}(w_0)\rangle\hookrightarrow P(w)$.
\item\label{p2.2} The restriction from $P$ to $T\langle-\mathfrak{l}(w_0)\rangle$ 
(the latter considered as a submodule of $P$ via the inclusion constructed in 
\eqref{p2.1}) induces an isomorphism,
\begin{displaymath}
\ttA=\End_{\ttA}(P)\cong \End_{\ttA}(T)=R(\ttA),
\end{displaymath}
of graded algebras.
\end{enumerate}
\end{proposition}

\begin{proof}
We start with the ungraded version, proved in \cite[Proposition~4]{FuKoMa}. 
Let $w=s_1\dots s_k$ be a reduced decomposition of $w$. Using the induction on
$\mathfrak{l}(w)$ one shows that applying $\theta_{s_k}\dots\theta_{s_1}$ to the
ungraded inclusion $T(w_0)\hookrightarrow P(e)$ induces the ungraded inclusion
$T(w_0w)\hookrightarrow P(w)$. This proves the ungraded analogue of \eqref{p2.1}
and the ungraded analogue of \eqref{p2.2} follows from \eqref{p2.1} using Enright's 
completion functor.

Obviously, the graded version of \eqref{p2.2} follows from the graded version of 
\eqref{p2.1} and the ungraded version of \eqref{p2.2}. Hence we are left to prove the
graded version of \eqref{p2.1}. We start with
\begin{equation}\label{eqeqp2}
T(w_0)\langle-\mathfrak{l}(w_0)\rangle\hookrightarrow P(e).
\end{equation}
Since $T(w_0)$ is the simple socle of $\Delta(x)$ for every $x\in\bfW$ we 
certainly have $T(w_0)\langle -k(x)\rangle\hookrightarrow \Delta(x)$ for some 
$k(x)$. In particular, $k(w_0)=0$ since $T(w_0)\cong \Delta(w_0)$. Using induction 
on $\mathfrak{l}(x)$, graded translation functors, and \cite[Theorem~3.6]{St1} one shows 
that $k(x)=\mathfrak{l}(w_0)-\mathfrak{l}(x)$, in particular, $k(e)=\mathfrak{l}(w_0)$.
Now \eqref{eqeqp2} follows from the observation that $P(e)\cong\Delta(e)$.

The rest follows again by induction on $\mathfrak{l}(w)$ applying the graded translation
$\theta_{s_k}\dots\theta_{s_1}$ to \eqref{eqeqp2} and using \cite[Theorem~3.6]{St1}.
\end{proof}

\begin{proposition}\label{p5}
$R(\ttA)$ is a positively graded algebra.
\end{proposition}

\begin{proof}
By definition (see \cite{Ri1}) every indecomposable summand of $T$ has both a
standard and a costandard (graded) filtration, that is a filtration, whose 
subquotients are standard and costandard modules respectively. Moreover, every 
indecomposable summand is self-dual (with respect to $\star$). Every 
(ungraded) morphism from $T(x)$ to $T(y)$, $x,y\in\bfW$, is a linear combination 
of morphisms, each of which is induced by a (unique up to a non-zero scalar) map 
from some subquotient of a standard filtration of $T(x)$ to some subquotient of 
a costandard filtration of $T(y)$. This means that the statement of the 
proposition follows from the following lemma.

\begin{lemma}\label{l6}
Let $x\in\bfW$. Then every subquotient of any standard filtration of 
$T(x)$, which is not isomorphic to $\Delta(x)$, has the form 
$\Delta(y)\langle l\rangle$ with $l>0$.
\end{lemma}

\begin{proof}
We prove this by a downward induction on $\mathfrak{l}(x)$ with the basis 
$x=w_0$ being obvious. Let now $x\in \bfW$ and $s$ be a simple reflection
such that $\mathfrak{l}(x)>\mathfrak{l}(xs)$. Consider the modules
$T(x)$ and $\theta_s T(x)$. Using \cite[Theorem~3.6]{St1} and the inductive 
assumption we obtain that 
\begin{enumerate}[(a)]
\item\label{l6.e1} every subquotient of any standard filtration of 
$\theta_s T(x)$ has the form $\Delta(y)\langle l\rangle$ with $l\geq 0$.
\end{enumerate}
The question is when we can get $l=0$? First of all, again by \cite[Theorem~3.6]{St1} 
we obtain that $\Delta(xs)$ occurs as a subquotient of any standard filtration of 
$\theta_s T(x)$. 

Fix now $y\neq xs$ such that $\Delta(y)$ occurs a subquotient of any standard 
filtration of $\theta_s T(x)$. Using \cite[Theorem~3.6]{St1} and the inductive assumption,
we get that every such occurrence comes from some occurrence of $\Delta(y)\langle 1\rangle$ 
as a subquotient of some standard filtration of $T(x)$. Moreover, $ys>y$. Let $m_y$ 
denote the multiplicity of $\Delta(y)\langle 1\rangle$ in $T(x)$.

First we claim that $m_y\leq \dim\ext^1_{\ttA}(L(y)\langle 1\rangle,L(x))$. Indeed,
let $f$ be a direct sum of all primitive idempotents of $R(\ttA)$, which correspond to 
$z\in\bfW$ such that $\mathfrak{l}(z)\geq \mathfrak{l}(x)$. By induction we can assume 
that the grading of $f R(\ttA) f$ is positive. Using Proposition~\ref{p2} and the graded 
contravariant Ringel duality functor $\hom_{\ttA}({}_-,\oplus_{l\in\Z}T\langle l\rangle)$, 
the statement reduces to the analogous statement for projective modules. But for projective 
modules over positively graded quasi-hereditary algebras the corresponding inequality 
$m_y\leq \dim\ext_{\ttA}^1(L(w_0x),L(w_0y)\langle -1\rangle)$ is obvious.

Now we recall the KL-combinatorics. From \cite[Proposition~5.2.3]{Ir}, the graded
Ringel self-duality, and the KL-Theorem it follows that we have the following 
decomposition in $\ttA\mathrm{-mod}$:
\begin{equation}\label{l6.e2}
\theta_s T(x)\cong T(xs)\oplus
\bigoplus_{y>x,ys>y}\left(\dim \Ext_{\ttA}^1(L(y),L(x))\right) T(y).
\end{equation}

Let us recall that all tilting modules are self-dual (even as graded modules).
This and \eqref{l6.e1} implies that \eqref{l6.e2} must be the case even in
$\ttA\mathrm{-grmod}$ (otherwise any shifted direct summand must come with an 
isomorphic direct summand shifted in the opposite way, which would imply that there
should occur some $\Delta(u)\langle t\rangle$, $t<0$, in any standard filtration
of $\theta_s T(x)$, contradicting \eqref{l6.e1}). Further, since the $s$-translation of
a standard filtration of $T(x)$ gives rise to a standard filtration of $\theta_s T(x)$ 
in a canonical way, it follows from \cite[Theorem~3.6]{St1} that in 
$\ttA\mathrm{-grmod}$ the multiplicity of $T(y)$ as a direct summand of 
$\theta_s T(x)$ can not exceed $m_y$. Since
\begin{displaymath}
m_y\leq \dim\ext^1_{\ttA}(L(y)\langle 1\rangle,L(x))\leq \dim\Ext_{\ttA}^1(L(y),L(x)),
\end{displaymath}
it follows that 
\begin{equation}\label{l6.e3}
m_y=\dim\ext^1_{\ttA}(L(y)\langle 1\rangle,L(x))=\dim\Ext_{\ttA}^1(L(y),L(x)),
\end{equation} 
which means that each occurrence of $\Delta(y)$ as a subquotient of a standard 
filtration of $\theta_s T(x)$ in fact comes from a direct summand of $\theta_s T(x)$, 
which is isomorphic to $T(y)$. This implies that every subquotient of any standard 
filtration of $T(xs)$, which is not isomorphic to $\Delta(xs)$, has the form 
$\Delta(z)\langle l\rangle$ with $l>0$, and completes the proof.
\end{proof}
\end{proof}

The proof of Proposition~\ref{p5} suggests the following immediate
corollary from \cite[Proposition~5.2.3]{Ir} and Proposition~\ref{p2}:

\begin{corollary}\label{s2c1}
Let $x\in\bfW$ and $s$ be a simple reflection. Then
\begin{displaymath}
\theta_s T(x)\cong
\begin{cases}
T(x)\langle 1\rangle\oplus T(x)\langle -1\rangle, & xs>x;\\
T(xs)\oplus
\bigoplus_{y>x,ys>y}\left(\dim \Ext_{\ttA}^1(L(y),L(x))\right) T(y), & xs<x.
\end{cases}
\end{displaymath}
\end{corollary}

Recall, that for a fixed graded module, $M$, over a graded algebra, $A$,
the complex $\mathcal{X}^{\bullet}$ is called {\em linear} provided that
$\mathcal{X}^{i}\in\mathrm{add}(M\langle i\rangle)$ for all $i\in\mathbb{Z}$.
Further, recall that the $\ttA$-module $N$ is said to have a {\em linear tilting 
(co)resolution} if the tilting (co)resolution of $N$ is a linear complex with 
respect to the graded lift of $T$ fixed above.

\begin{proposition}\label{s2p1}
The module $\Delta$ admits a (finite) linear tilting coresolution and the module 
$\nabla$  admits a (finite) linear tilting resolution.
\end{proposition}

\begin{proof}
By duality it is enough to prove the statement for $\Delta$. It
is further enough to show that $\Delta(x)$ admits a linear tilting coresolution
for every $x\in \bfW$. We show this by induction on $\mathfrak{l}(w_0)-
\mathfrak{l}(x)$. The basis of the induction is obvious since  $\Delta(w_0)$ is 
a tilting module. Let $x\in\bfW$ and $s$ be a simple reflection such that
$xs<x$. By induction we can assume that $\Delta(x)$ has the linear tilting 
coresolution $\mathcal{T}^{\bullet}(\Delta(x))$. Applying $\theta_s$ we obtain
a tilting coresolution of $\theta_s\Delta(x)$ and, since all tilting modules
have costandard filtrations, the adjunction induces a morphism
of complexes, $\varphi^{\bullet}:\theta_s\mathcal{T}^{\bullet}(\Delta(x))\tto 
\mathcal{T}^{\bullet}(\Delta(x))\langle 1\rangle$, surjective on every component.
Let $\mathcal{Q}^{\bullet}$ denote the cone of $\varphi^{\bullet}$ shifted by
$[-1]$. Then $\mathcal{Q}^{\bullet}$ is a tilting coresolution of $\Delta(xs)$. 

We claim that, taking away all trivial direct summands of $\mathcal{Q}^{\bullet}$ 
(that is direct summands of the form $\dots\to 0\to T(z)\overset{\mathrm{id}}{\to} 
T(z)\to 0\to\dots$), $\mathcal{Q}^{\bullet}$ reduces to a linear complex, $\overline{\mathcal{Q}}^{\bullet}$, 
of tilting modules. Indeed, let us fix $l\geq 0$ and let $\mathcal{T}^{l}(\Delta(x))\cong
\oplus_{w\in\bfW} m_w  T(w)\langle l\rangle$. Then, by Corollary~\ref{s2c1}, we have
\begin{displaymath}
\theta_s\mathcal{T}^{l}(\Delta(x))\cong
\oplus_{ws>w} m_w  T(w)\langle l+1\rangle\bigoplus \oplus_{ws>w} m_w  
T(w)\langle l-1\rangle \bigoplus X,
\end{displaymath}
where $X\in \mathrm{add}(T\langle l\rangle)$. It is easy to see that the morphism 
$\varphi^{\bullet}$ induces an isomorphism  between the corresponding direct summands
$T(w)\langle l+1\rangle$ in $\theta_s\mathcal{T}^{l}(\Delta(x))$ and
$\mathcal{T}^{l}(\Delta(x))\langle 1\rangle$, implying that
$\overline{\mathcal{Q}}^{l}\in \mathrm{add}(T\langle l\rangle\oplus T\langle l-1\rangle)$.

On the other hand, let $f^{\bullet}$ be the differential in 
$\mathcal{T}^{\bullet}(\Delta(x))$. Then
each occurrence of $T(w)\langle l\rangle$ in $\mathcal{T}^{l}(\Delta(x))$ comes from some 
standard subquotient $\Delta(w)\langle l\rangle$ in the cokernel of $f^{l-2}$ and thus in
$\mathcal{T}^{l-1}(\Delta(x))$. Hence, by Corollary~\ref{s2c1}, in the case $ws>w$  it 
gives rise to an occurrence of $T(w)\langle l-1\rangle$ in 
$\theta_s\mathcal{T}^{l-1}(\Delta(x))$. Further, we obtain that $\theta_s f^{l}$ induces 
a non-zero map and hence an isomorphism between the corresponding summands 
$T(w)\langle l-1\rangle$ in $\theta_s\mathcal{T}^{l-1}(\Delta(x))$ and
in $\theta_s\mathcal{T}^{l}(\Delta(x))$. This splits away and thus we obtain that 
$\overline{\mathcal{Q}}^{l}\in \mathrm{add}(T\langle l\rangle)$, completing the proof.
\end{proof}

\begin{corollary}\label{s2p1c1}
The module $\Delta$ admits a (finite) linear projective resolution and the module 
$\nabla$  admits a (finite) linear injective coresolution.
\end{corollary}

\begin{proof}
This follows from Proposition~\ref{s2p1} and the graded Ringel self-duality of $\ttA$.
\end{proof}

\begin{proof}[Proof of Theorem~\ref{t1}.]
From Proposition~\ref{p5} we already know that the induced grading on
$R(\ttA)$ is positive. Hence the natural grading on $\ttA$ is positive by
Proposition~\ref{p2}\eqref{p2.2}. By Corollary~\ref{s2p1c1}, every standard 
$\ttA$-module admits a linear projective resolution, and every costandard 
$\ttA$-module admits a linear injective coresolution. This means that
$\ttA$ is a standard Koszul quasi-hereditary algebra with respect to the
underlying (natural) grading. The theorem is proved.
\end{proof}

\section{The category of linear complexes of tilting modules}\label{s3}

\subsection{Quadratic dual via linear complexes}\label{s3.1}

Let $A=\oplus_{i\geq 0}A_i$ be a basic positively graded algebra over some field 
$\Bbbk$ with finite-dimensional graded components. Denote by $A\mathrm{-fgmod}$ the 
category of all graded $A$-modules with finite-dimensional graded components. Let us 
fix the natural gradings on simple and projective $A$-modules, which is induced
from the positive grading on $A$. Let $\mathcal{P}(A)$ denote the
category, whose objects are all linear complexes of finitely generated
projective $A$-modules, that is $\mathcal{X}^{\bullet}\in \mathcal{P}(A)$ 
if and only if $\mathcal{X}^{l}\in \mathrm{add}\left({}_A A\langle l\rangle\right)$
for all $l\in\Z$, and morphism are all possible morphisms of complexes of graded
modules. For an $A$-module, $M$, we denote by $M^{\bullet}$ the complex satisfying 
$M^{0}=M$ and $M^l=0$ for all $l\neq 0$. For a $\Bbbk$-vector space, $V$, we set
$V^*=\mathrm{Hom}_{\Bbbk}(V,\Bbbk)$. For two $\Bbbk$-vector space, $V$ and $W$,
and for $f\in \mathrm{Hom}_{\Bbbk}(V,W)$ we denote by $f^*$ the corresponding
dual map from $\mathrm{Hom}_{\Bbbk}(W^*,V^*)$. For finite-dimensional $V$ and $W$
we have an obvious canonical isomorphism, $(V\otimes_{\Bbbk}W)^*\cong 
V^*\otimes_{\Bbbk}W^*$.

Denote further by $A^!$ the {\em quadratic dual} of $A$, defined as follows:
$A^!=A_0[A_1^*]/\left(\mu^*(A_2^*)\right)$, where $\mu:A_1\otimes_{A_0}A_1\to A_2$ 
is the multiplication map. Note that $A^!$ is quadratic by definition. 
If $A$ is quadratic we have $A\cong \left(A^!\right)^!$ canonically. Recall the
following statement (see \cite[Theorem~2.4]{MaSa} or \cite[Theorem~7]{MaOv}):

\begin{theorem}\label{tdual}
There is an equivalence, $F: A^!\mathrm{-fgmod}\cong \mathcal{P}(A)$.
\end{theorem}

Let $e\in A_0$ be an idempotent and let $\mathcal{P}^e(A)$ denote the full subcategory
of $\mathcal{P}(A)$ which consists of all $\mathcal{X}^{\bullet}\in \mathcal{P}(A)$ 
such that  $\mathcal{X}^{l}\in \mathrm{add}\left({}_A Ae\langle l\rangle\right)$.
Note that there is a canonical bijection between the idempotents in $A$ and $A^!$
(since both $A$ and $A^!$ share the same semi-simple part $A_0$ by definition).
In particular, we can define $B_e=A^!/(A^! (1-e) A^!)$. In what follows we will need 
the following easy corollary from Theorem~\ref{tdual}:

\begin{proposition}\label{pdual}
$F$ induces the equivalence $F^e:(eAe)^!\mathrm{-fgmod}\cong
\mathcal{P}^e(A)$. In particular, $(eAe)^!\cong B_e$ as graded algebras with respect 
to the induced gradings.
\end{proposition}

\begin{proof}
Every projective object from $\mathcal{P}^e(A)$ is by definition
the maximal quotient of the corresponding projective object from 
$\mathcal{P}(A)$, which contains only simple subquotients of the form 
$\mathrm{add}\left({}_A Ae^{\bullet}\langle l\rangle\right)$. This
and \cite[Section~5]{Au} imply the first statement, and the second statement 
follows from the first one.
\end{proof}

\subsection{Application to the regular block of $\ccO$}\label{s3.2}

Let $\cT(\ttA)$ denote the category, whose objects are all bounded complexes 
$\mathcal{X}^{\bullet}$ of graded $\ttA$-modules satisfying
$\mathcal{X}^{l}\in\mathrm{add}\, (T\langle l\rangle)$ for all $l\in\Z$
(such complexes are called {\em linear complexes of tilting modules}),
and morphisms are all usual morphisms of complexes of graded modules, 
see \cite{MaOv}. We also denote by $\cP(\ttA)$ the category, whose objects 
are all bounded complexes $\mathcal{X}^{\bullet}$ of graded $\ttA$-modules 
satisfying $\mathcal{X}^{l}\in\mathrm{add}\,(P\langle l\rangle)$ for all $l\in\Z$,
and morphisms are all usual morphisms of complexes of graded modules
(such complexes are called {\em linear complexes of projective modules}), 
see \cite{MaOv}. 

\begin{theorem}\label{t21}
\begin{enumerate}[(i)]
\item\label{t21.1} There is an equivalence, 
$\mathcal{F}:\ttA\mathrm{-fgmod}\cong\cT(\ttA)$, which
sends $L(w)$ to $T(w_0 w^{-1} w_0)^{\bullet}$ for every $w\in\bfW$.
\item\label{t21.2} There is an equivalence, 
$\mathcal{G}:\ttA\mathrm{-fgmod}\cong\cP(\ttA)$, which
sends $L(w)$ to $P(w^{-1} w_0)^{\bullet}$ for every $w\in\bfW$.
\end{enumerate}
\end{theorem}

\begin{proof}
By Theorem~\ref{tdual} there is an equivalence, $\ttA^!\mathrm{-fgmod}\cong\cP(\ttA)$. 
Applying the graded version of the Ringel duality from \cite{So2} gives an equivalence, 
$R(\ttA)^!\mathrm{-fgmod}\cong\cT(\ttA)$. Because of Proposition~\ref{p2} we even 
obtain an equivalence, $\ttA^!\mathrm{-fgmod}\cong\cT(\ttA)$. However, $\ttA$ is Koszul 
and even Koszul self-dual by \cite{So1} and the natural grading on $\ttA$ is the Koszul 
one by Theorem~\ref{t1}. Thus $\ttA^!\cong \ttA$ by \cite[2.9]{BeGiSo}, which proves 
the existence of both $\mathcal{F}$ and $\mathcal{G}$. The correspondence on simple 
objects follows from Proposition~\ref{p2}\eqref{p2.1} and \cite{So1}.
\end{proof}

\section{Singular blocks of $\ccO$ and their Koszulity}\label{s4}

Let $\bfG\subset\bfW$ be a parabolic subgroup and 
$w_0^{\bfG}$ be the longest element in $\bfG$. Denote by $\bfW_{\bfG}$ the 
set of the longest coset representatives in $\bfW/\bfG$. Let $\lambda\in\fh^*$ be a 
dominant integral (singular) weight with stabilizer $\bfG$. Denote by $\cO_{\lambda}$ 
the singular block of $\cO$, which corresponds to $\lambda$, and by $\ttA_{\bfG}$ 
the corresponding basic associative algebra. The bijection between the simple
$\ttA$-modules and the elements of $\bfW$, described in the introduction,
induces a bijection between the isomorphism classes of simple $\ttA_{\bfG}$-modules 
and the cosets from $\bfW/\bfG$, and, also, the elements of $\bfW_{\bfG}$. 
For $\ttA_{\bfG}$-modules we will use notation $L_{\bfG}(w)$, 
$w\in \bfW_{\bfG}$, etc. Soergel's combinatorics from \cite{So1} equips 
$\ttA_{\bfG}$ with a natural grading in the following way. The algebra
$\End_{\ttA_{\bfG}}\left(P_{\bfG}(w_0)\right)$ is the subalgebra $\bfC_{\bfG}$ of 
$\bfG$-invariants in $\bfC$, in particular, is graded. Moreover, the module
$P_{\bfG}(w_0)_{\bfC_{\bfG}}$ is a graded module,. This induces a grading on
$\ttA_{\bfG}\cong \End_{\bfC_{\bfG}}\left(P_{\bfG}(w_0)_{\bfC_{\bfG}}\right)$,
which we call {\em natural}. Let $\theta_{\bfG}^{on}:\cO_0\to \cO_{\lambda}$ and 
$\theta_{\bfG}^{out}:\cO_{\lambda}\to \cO_0$ denote the functors of translations onto 
and out of the $\bfG$-wall respectively. These functors are left and right adjoint to 
each other. Set $\theta_{\bfG}=\theta_{\bfG}^{out}\theta_{\bfG}^{on}$. 
In the same way as it is done in Section~\ref{s2} we fix the graded lifts of 
simple, projective, injective, standard, costandard and tilting $\ttA_{\bfG}$-modules
and will use for them analogous notation (for example $L_{\bfG}(w)$ etc.).
The main result of the present section is the following statement.

\begin{theorem}\label{s2.2.t2}
The algebra $\ttA_{\bfG}$ is a standard Koszul quasi-hereditary algebra, 
and the natural grading on $\ttA_{\bfG}$ is the Koszul one.
\end{theorem}

Again, to prove this theorem we will need some preparation.

\begin{lemma}\label{s2.2.l201}
The natural grading on $\ttA_{\bfG}$ is positive.
\end{lemma}

\begin{proof}
Let $x\in \bfW_{\bfG}$ and $y\in \bfW$. Then the adjointness of
$\theta_{\bfG}^{out}$ and  $\theta_{\bfG}^{on}$ gives
\begin{equation}\label{eqs2.2.l201}
\Hom_{\ttA}\left(\theta_{\bfG}^{out}P_{\bfG}(x),L(y)\right)\cong
\Hom_{\ttA_{\bfG}}\left(P_{\bfG}(x),\theta_{\bfG}^{on}L(y)\right).
\end{equation}
However, since $\theta_{\bfG}^{on}$ sends simple $\ttA$-modules to 
simple $\ttA_{\bfG}$-modules or zero, it follows that for a fixed $y$
the right-hand side of \eqref{eqs2.2.l201} is either $0$ for all $x$, or
is $0$ for all $x$ except one, for which it is equal to $\mathbb{C}$.
This implies that the projective module $\theta_{\bfG}^{out}P_{\bfG}(x)$ 
is in fact indecomposable. Using \cite[Theorem~8.2]{St1} it is easy to 
see that $\theta_{\bfG}^{out}$ admits a natural graded lift compatible with 
the grading on both $\cO_{\lambda}$ and $\cO_0$. In particular, it follows 
that the algebra $\ttA_{\bfG}$ is a graded subalgebra of $\ttA$ and hence 
Theorem~\ref{t1} implies that $\ttA_{\bfG}$ is positively graded. 
\end{proof}

\begin{proposition}\label{s2.2.l2}
The module $\theta_{\bfG}^{out}\Delta_{\bfG}(w)\langle \mathfrak{l}(w_0^{\bfG})\rangle$, 
$w\in \bfW_{\bfG}$, admits a linear tilting coresolution in $\ttA\mathrm{-grmod}$.
\end{proposition}

\begin{proof}
Throughout the proof we fix $w\in\bfW_{\bfG}$. 
For $i=0,1,\dots,\mathfrak{l}(w_0^{\bfG})$ let 
\begin{gather*}
W_i^w=\{x\in w_0\bfG w^{-1}w_0\,:\, 
\mathfrak{l}(x)=\mathfrak{l}(w^{-1})-\mathfrak{l}(w_0^{\bfG})+i\},\\
\overline{W}_i^w=\{x\in w\bfG\,:\, \mathfrak{l}(x)=\mathfrak{l}(w)-
\mathfrak{l}(w_0^{\bfG})+i\}
\end{gather*}
(in particular, $W_0^w=\{w_0w_0^{\bfG}w^{-1}w_0\}$ and $\overline{W}_0^w=\{ww_0^{\bfG}\}$).
For $w\in\bfW$ we denote by $\mathcal{T}^{\bullet}(\Delta(w))$ the 
linear tilting resolution of $\Delta(w)$. 

For $i=0,1,\dots,\mathfrak{l}(w_0^{\bfG})$ let $X_i$ denote the (non-graded) trace in
$\theta_{\bfG}^{out}\Delta_{\bfG}(w)$ of all $P(x)$ such that there is $y\in \overline{W}_i^w$ satisfying $y\geq x$. We consider the trace for non-graded modules, however, 
since both $\theta_{\bfG}^{out}\Delta_{\bfG}(w)$ and $P(x)$ are graded, the trace
itself will be a graded submodule of $\theta_{\bfG}^{out}\Delta_{\bfG}(w)$,
see \cite[Lemma~4]{MaOv}. Set $Y_i=\left(\theta_{\bfG}^{out}\Delta_{\bfG}(w)\right)/X_i$. 
In particular, we have $X_0\cong \Delta(ww_0^{\bfG})$ and $Y_{\mathfrak{l}(w_0^{\bfG})-1}\cong
\Delta(w)$.

\begin{lemma}\label{help1}
For all $x\in w\bfG$ and $k>0$ we have the following equality:
$\Ext_{\ttA}^{k}\left(\theta_{\bfG}^{out}\Delta_{\bfG}(w),\Delta(x)\right)=0$.
\end{lemma}

\begin{proof}
The functors $\theta_{\bfG}^{out}$, $\theta_{\bfG}^{on}$ are exact, left and right 
adjoint to each other, and send projectives to projectives. Hence
\begin{multline*}
\Ext_{\ttA}^{k}\left(\theta_{\bfG}^{out}\Delta_{\bfG}(w),\Delta(x)\right)=
\Ext_{\ttA}^{k}\left(\Delta_{\bfG}(w),\theta_{\bfG}^{on}\Delta(x)\right)=\\=
\Ext_{\ttA}^{k}\left(\Delta_{\bfG}(w),\Delta_{\bfG}(w)\right)=0
\end{multline*}
since Verma modules do not have self-extensions in $\cO$.
\end{proof}

\begin{lemma}\label{help2}
$\Ext_{\ttA}^1\left(Y_{i},\Delta(x)\right)=\mathbb{C}$ for 
$i\in\{0,1,\dots,\mathfrak{l}(w_0^{\bfG})-1\}$ and
$x\in \overline{W}_{i}^w$.
\end{lemma}

\begin{proof}
Applying $\Hom_{\ttA}\left({}_-,\Delta(x)\right)$ to
the short exact sequence 
\begin{displaymath}
X_i\hookrightarrow \theta_{\bfG}^{out}\Delta_{\bfG}(w)\tto Y_i
\end{displaymath}
we obtain the following fragment in the long exact sequence:
\begin{multline}\label{helpeq}
0\to
\Hom_{\ttA}\left(Y_i,\Delta(x)\right)\overset{\sim}{\to}
\Hom_{\ttA}\left(\theta_{\bfG}^{out}\Delta_{\bfG}(w),\Delta(x)\right)\to\\ \to
\Hom_{\ttA}\left(X_i,\Delta(x)\right)\to
\Ext_{\ttA}^1\left(Y_{i},\Delta(x)\right)\to
\Ext_{\ttA}^1\left(\theta_{\bfG}^{out}\Delta_{\bfG}(w),\Delta(x)\right)\to
\end{multline}
Since $Y_i$ surjects, by the definition, onto $\Delta(w)$, which, in turn,
is a submodule of $\Delta(x)$ by \cite[7.7.7]{Di}, we get 
$\Hom_{\ttA}\left(Y_i,\Delta(x)\right)\neq 0$. Using the adjointness of
$\theta_{\bfG}^{out}$ and $\theta_{\bfG}^{on}$, we obtain
\begin{multline*}
\Hom_{\ttA}\left(\theta_{\bfG}^{out}\Delta_{\bfG}(w),\Delta(x)\right)\cong
\Hom_{\ttA}\left(\Delta_{\bfG}(w),\theta_{\bfG}^{on}\Delta(x)\right)\cong\\\cong
\Hom_{\ttA}\left(\Delta_{\bfG}(w),\Delta_{\bfG}(w)\right)\cong\mathbb{C},
\end{multline*}
which establishes the isomorphism indicated in \eqref{helpeq}. 
Further, the only standard subquotient of $X_i$, which maps to
$\Delta(x)$ is $\Delta(x)$ itself, which implies 
$\Hom_{\ttA}\left(X_i,\Delta(x)\right)\cong\mathbb{C}$. Finally,
$\Ext_{\ttA}^1\left(\theta_{\bfG}^{out}\Delta_{\bfG}(w),\Delta(x)\right)=0$
by Lemma~\ref{help1}. The exactness of \eqref{helpeq} now gives the
necessary statement.
\end{proof}

Lemma~\ref{help2} implies that the modules $Y_i$ can be constructed recursively
starting from $i=\mathfrak{l}(w_0^{\bfG})$ and descending to
$i=0$. On every step of this recursion one uses the universal extension procedure
as follows: we start with $Y_{\mathfrak{l}(w_0^{\bfG})}\cong\Delta(w)$ and proceed,
involving the module $\bigoplus_{x\in \overline{W}_{i}^w}\Delta(x)$ on the step 
$\mathfrak{l}(w_0^{\bfG})-i$. Note that non-isomorphic direct summands of the 
module  $\bigoplus_{x\in \overline{W}_{i}^w}\Delta(x)$ 
do not have self-extensions, see for example \cite[Section~9]{RoCa}.

Let $\mathfrak{a}$ be the semi-simple Lie subalgebra of $\mathfrak{g}$, 
associated with $\bfG$. Applying the parabolic induction to the classical 
BGG-resolution of a finite-dimensional $\mathfrak{a}$-module gives the following 
resolution,
\begin{equation}\label{s2.2.l2.e101}
0\to \bigoplus_{x\in W_{\mathfrak{l}(w_0^{\bfG})}^w}\Delta(x)\to\dots\to 
\bigoplus_{x\in W_1^w}\Delta(x)\to\Delta(w_0w_0^{\bfG}w^{-1}w_0)\to 0,
\end{equation}
of the generalized Verma module, which is the cokernel of the last non-zero map
in \eqref{s2.2.l2.e101}. By \cite[Theorem~8]{MaOv}, the sequence 
\eqref{s2.2.l2.e101} induces, via the equivalence $\mathcal{F}$ from 
Theorem~\ref{t21}, the following complex of elements from $\cT(\ttA)$:
\begin{equation}\label{s2.2.l2.e1}
0\to \bigoplus_{x\in \overline{W}_{\mathfrak{l}(w_0^{\bfG})}^w}
\mathcal{T}^{\bullet}(\Delta(x))\to\dots\to 
\bigoplus_{x\in \overline{W}_1^w}\mathcal{T}^{\bullet}(\Delta(x))
\to\mathcal{T}^{\bullet}(\Delta(ww_0^{\bfG}))\to 0.
\end{equation}
The sequence \eqref{s2.2.l2.e1} is exact in all terms except the last one.
This allows us to take inductively the cone of all morphisms in
\eqref{s2.2.l2.e1} and, moreover, implies that on every step the complex
of tilting modules we obtain is isomorphic to a linear complex of tilting 
modules. It is easy to see that the homology of the complex, obtained on every
step, is concentrated in a single position (in particular, in position
$0$ on the last step). Let us denote this homology by $Z_i$,
$i=\mathfrak{l}(w_0^{\bfG}),\mathfrak{l}(w_0^{\bfG})-1,\dots$.

Let us show, by induction on $i$, that $Z_i\cong Y_{i-1}$.
This is obvious for $i=\mathfrak{l}(w_0^{\bfG})$. Note that the restriction 
of all differentials in the BGG resolution to Verma modules are non-zero 
(see for example \cite[Section~10]{RoCa}). Let us now assume that
$Z_i\cong Y_{i-1}$. In particular, we have
$\Ext_{\ttA}^1\left(Z_{i},\Delta(x)\right)=\mathbb{C}$ for every
$x\in \overline{W}_{i-1}^w$ by Lemma~\ref{help2}. The cone construction of
$Z_{i-1}$ now implies that $Z_{i-1}$ is obtained from $Z_i$ via the universal
extension with $\bigoplus_{x\in W_{i-1}^w}\Delta(x)$. By the universality of
this extension we obtain that $Z_{i-1}$ must be isomorphic to $Y_{i-2}$.

In particular, $Z_{0}\cong \theta_{\bfG}^{out}\Delta_{\bfG}(w)
\langle\mathfrak{l}(w_0^{\bfG})\rangle$ and this module has the linear tilting 
coresolution as prove above.
\end{proof}

\begin{lemma}\label{nnlemn}
\begin{enumerate}[(i)]
\item\label{nnlemn.1} Let $w\in\bfW_{\bfG}$. Then the (ungraded)
trace $Tr_{\bfG}(w)$ of $P_{\bfG}(w_0)$ in $P_{\bfG}(w)$ is an indecomposable 
tilting $\ttA_{\bfG}$-module. Moreover, for $w\neq w'$, $w'\in\bfW_{\bfG}$, we have  
$Tr_{\bfG}(w)\not\cong Tr_{\bfG}(w')$.
\item\label{nnlemn.2} Restriction from $P_{\bfG}$ to the trace of
$P_{\bfG}(w_0)$ in $P_{\bfG}$ establishes the Ringel self-duality of $\ttA_{\bfG}$.
\end{enumerate}
\end{lemma}

\begin{proof}
From Proposition~\ref{p2} it follows that $T(w_0w)$ is the trace of $P(w_0)$ in
$P(w)$. $\theta_{\bfG}^{on}$ sends $L(w_0)$ to $L_{\bfG}(w_0)$, which implies that 
the trace of $P(w_0)$ in $P(w)$ is mapped to the trace of $\theta_{\bfG}^{on} P(w_0)$ 
in $\theta_{\bfG}^{on} P(w)$. This and the fact that $\theta_{\bfG}^{on}$ sends 
tilting modules to tilting modules implies that $Tr_{\bfG}(w)$ is a tilting module.
On the other hand, from \cite[3.1]{KoSlXi} we have that every $P_{\bfG}(w)$ has a 
two-step copresentation by modules from $\mathrm{add}\left(P_{\bfG}(w_0)\right)$. 
This implies that the Auslander $P_{\bfG}(w_0)$-coapproximation of $Tr_{\bfG}(w)$, 
see \cite[Section~5]{Au}, is isomorphic to $P_{\bfG}(w)$. Using the same arguments 
as in \cite[Theorem~4.1]{FuKoMa} we obtain that the restriction defines an 
isomorphism between the algebras $\mathrm{End}_{\ttA_{\bfG}}\left(P_{\bfG}\right)$ and
$\mathrm{End}_{\ttA_{\bfG}}\left(\oplus_{w\in \bfW_{\bfG}}Tr_{\bfG}(w)\right)$. 
Both \eqref{nnlemn.1} and \eqref{nnlemn.2} follow.
\end{proof}

\begin{proof}[Proof of Theorem~\ref{s2.2.t2}]
Observe that $\theta_{\bfG}^{out}$ sends tilting modules from $\cO_{\lambda}$
to tilting modules from $\cO_0$, and projective modules from $\cO_{\lambda}$
to projective modules from $\cO_0$. Adjointness with $\theta_{\bfG}^{on}$ 
implies that $\theta_{\bfG}^{out}$ sends indecomposable projective modules
to indecomposable projective modules. This and Lemma~\ref{nnlemn} gives that 
$\theta_{\bfG}^{out}$ even sends indecomposable tilting modules to 
indecomposable tilting modules. In particular, for every $w\in \bfW_{\bfG}$ the 
minimal tilting coresolution of $\Delta_{\bfG}(w)$ is sent to a minimal tilting 
coresolution of $\theta_{\bfG}^{out}\Delta_{\bfG}(w)$. The last one is however
linear (up to a shift of grading) by Proposition~\ref{s2.2.l2}. 
Since $\theta_{\bfG}^{out}$ is compatible with
the gradings on $\cO_{\lambda}$ and $\cO_0$, we obtain that the original tilting
coresolution of $\Delta_{\bfG}(w)$ was linear. From \cite[Theorem~6]{MaOv}
it now follows that $\ttA_{\bfG}$ is standard Koszul with respect to the
natural grading.
\end{proof}

\section{The parabolic category of Rocha-Caridi \\ and its Koszulity}\label{s5}

Let now $\bfG$ and $\bfH$ be two parabolic subgroups of $\bfW$, and 
$\lambda$ and $\cO_{\lambda}$ be as in Section~\ref{s4}. Let further $\bfW_{\bfG}^{\bfH}$ 
denote the set of all $w\in\bfW$ which are the longest coset representatives in 
$\bfW/\bfG$ and the shortest coset representatives in $\bfH\backslash\bfW$ at the same 
time. Let $f_{\bfG}^{\bfH}$ be the sum of all primitive idempotents of $\ttA_{\bfG}$, which 
correspond to $w\in \bfW_{\bfG}^{\bfH}$ and set
$\ttA_{\bfG}^{\bfH}=\ttA_{\bfG}/\left(\ttA_{\bfG}(1-f_{\bfG}^{\bfH})\ttA_{\bfG}\right)$.
Then the algebra $\ttA_{\bfG}^{\bfH}$ is the associative algebra of  
Rocha-Caridi's parabolic subcategory of $\cO_{\lambda}$, associated with
$\bfH$. This is the full subcategory o $\cO_{\lambda}$, which consists of all 
modules, which are locally finite with respect to the parabolic subalgebra 
$\mathfrak{p}$ of $\fg$, associated with $\bfH$. Let $\mathfrak{a}$ be the
semi-simple part of $\mathfrak{p}$. The natural grading on 
$\ttA_{\bfG}$ induces a grading on $\ttA_{\bfG}^{\bfH}$ in an obvious way. We call 
the later grading again {\em natural}. The $\ttA_{\bfG}^{\bfH}$-modules will be denoted
$L_{\bfG}^{\bfH}(w)$ etc. The main result of this section is the 
following statement.

\begin{theorem}\label{s5.t1}
The algebra $\ttA_{\bfG}^{\bfH}$ is a standard Koszul quasi-hereditary algebra
and the natural grading on $\ttA_{\bfG}^{\bfH}$ is the Koszul one.
\end{theorem}

\begin{proof}
It is clear that the natural grading on $\ttA_{\bfG}^{\bfH}$ is positive
because the natural grading on $\ttA_{\bfG}$ is positive. Hence, to complete
the proof it would suffice to show that standard $\ttA_{\bfG}^{\bfH}$-modules admit
linear tilting coresolution and that costandard $\ttA_{\bfG}^{\bfH}$-modules admit
linear tilting resolution and apply \cite[Theorem~6]{MaOv}.

To show this we first forget about the grading for a moment.
Let $Z_{\bfH}:\ttA_{\bfG}\mathrm{-mod}\to \ttA_{\bfG}^{\bfH}\mathrm{-mod}$ 
denote Zuckerman's functor of taking the maximal subquotient, which belongs
to $\ttA_{\bfG}^{\bfH}\mathrm{-mod}$. Let $x\in \bfW_{\bfG}^{\bfH}$ and 
\begin{equation}\label{s4p1.e1}
0\to P_k\to\dots\to P_0\to\Delta_{\bfG}(x)\to 0
\end{equation}
be a projective resolution of $\Delta_{\bfG}(x)\in\ttA_{\bfG}\mathrm{-mod}$. 
Let us look at \eqref{s4p1.e1} inside $\cO_{\lambda}$ instead of 
$\ttA_{\bfG}\mathrm{-mod}$ keeping the notation. Since every $x\in \bfW_{\bfG}^{\bfH}$
corresponds to an $\mathfrak{a}$-dominant weight by definition, the corresponding
Verma  module $\Delta_{\bfG}(x)$ is obtained using the parabolic 
induction from some dominant Verma $\mathfrak{a}$-modules. However, the dominant
Verma module is projective in the category $\cO$ for the algebra $\mathfrak{a}$
(we denote this category by $\cO(\mathfrak{a})$). Since $U(\mathfrak{g})$, considered as
an $\mathfrak{a}$-module under the adjoint action, is a direct sum of finite-dimensional 
modules, the parabolic induction reduces (on the level of $\mathfrak{a}$-modules) to 
the tensoring with finite-dimensional modules. Hence we obtain that 
$\Delta_{\bfG}(x)$ (in particular, $\Delta_{\bfG}(w_0^{\bfG})$) is a (possibly 
infinite) direct sum of projective modules from $\cO(\mathfrak{a})$. Since all 
projective modules in $\cO_{\lambda}$ have Verma flag, it follows that
all  projective modules in $\cO_{\lambda}$ are (possibly infinite) direct sums of 
projective modules from $\cO(\mathfrak{a})$ as well. 

This means that \eqref{s4p1.e1}, considered as a sequence of $\mathfrak{a}$-modules, 
is a direct sum of trivial sequences of the form 
$\dots\to 0\to X\overset{\mathrm{id}}{\to} X\to 0\to\dots$. 
Since $Z_{\bfH}$ commutes with the restriction to $\mathfrak{a}$ (by definition), 
we obtain that the sequence
\begin{equation}\label{s4p1.e2}
0\to Z_{\bfH} P_k\to\dots\to Z_{\bfH} P_0\to Z_{\bfH} \Delta_{\bfG}(x)\to 0
\end{equation}
is exact. Since $Z_{\bfH}$ sends projective $\ttA_{\bfG}$-modules to projective
$\ttA_{\bfG}^{\bfH}$-modules or zero, \eqref{s4p1.e2} is in fact a projective 
resolution of the standard $\ttA_{\bfG}^{\bfH}$-module $Z_{\bfG} \Delta_{\bfG}(x)$. 

Now let us go back to the graded algebras.
From Theorem~\ref{s2.2.t2} we know that $\ttA_{\bfG}$ is standard Koszul, 
in particular, all standard $\ttA_{\bfG}$-modules admit linear projective 
resolutions. This means that \eqref{s4p1.e1} can be chosen linear.
Since $Z_{\bfH}$ obviously respects the grading, we get that \eqref{s4p1.e2} is
linear as well. Hence all standard $\ttA_{\bfG}^{\bfH}$-modules admit linear
projective resolutions. Applying $\star$ we get the dual statement for the
injective modules. This shows that $\ttA_{\bfG}^{\bfH}$ is standard Koszul, 
in particular, Koszul, by \cite[Theorem~1]{ADL}, and completes the proof.
\end{proof}

\begin{remark}
We remark that for the algebra $\ttA_{\{e\}}^{\bfH}$ the proof of
Theorem~\ref{s5.t1} does not use Section~\ref{s4} and hence does not use
Soergel's Koszul self-duality result.
\end{remark}

\section{A new parabolic generalization of the category $\ccO$}\label{s6}

In this section we will develop one auxiliary tool, which we will later use in
Section~\ref{s7}. This is a new parabolic generalization of the category $\ccO$,
which is not a highest weight category in general, but rather corresponds to a 
properly stratified algebra. Roughly speaking, it the the translation of a singular
block of $\cO$ (or some parabolic category) out of the wall.

Let $\bfG$ and $\bfH$ be two parabolic subgroups of $\bfW$. Let further $f_{\bfG}^{\bfH}$ 
be the sum of all primitive idempotents of $\ttA_{\{e\}}^{\bfH}$, which correspond 
to $w\in \bfW_{\bfG}^{\bfH}$, and set
$\ttB_{\bfG}^{\bfH}=f_{\bfG}^{\bfH}\ttA_{\{e\}}^{\bfH}f_{\bfG}^{\bfH}$.
Let $\bfC(\bfG)$ denote the coinvariant algebra of $\bfG$, which we consider as a graded
algebra with respect to the grading for which the generators have degree two.

\begin{theorem}\label{t11}
\begin{enumerate}[(i)]
\item\label{t11.1} There is an isomorphism of graded algebras,
\begin{displaymath}
\ttB^{\bfH}_{\bfG}\cong \ttA^{\bfH}_{\bfG}\otimes_{\C} \bfC(\bfG);
\end{displaymath}
\item\label{t11.2} the algebra $\ttB^{\bfH}_{\bfG}$ is properly stratified 
(in general not quasi-hereditary) and has a duality;
\item\label{t11.3} $\ttB^{\bfH}_{\bfG}\mathrm{-mod}$ is equivalent to the 
full subcategory $\cX$ of $\ttA_{\{e\}}^{\bfH}\mathrm{-mod}$, which consists of all 
$X$ admitting a two-step projective presentation by modules from 
$\mathrm{add}\,\left(\ttA_{\{e\}}^{\bfH} f_{\bfG}^{\bfH}\right)$;
\item\label{t11.4} $\theta_{\bfG}^{out}$ sends standard $\ttA^{\bfH}_{\bfG}$-modules 
to standard $\ttB^{\bfH}_{\bfG}$-modules;
\item\label{t11.5} as objects in $\cX$, the proper costandard 
$\ttB^{\bfH}_{\bfG}$-modules are $\Delta^{\bfH}_{\{e\}}(ww_{0}^{\bfG})^{\star}$, where 
$w\in \bfW_{\{e\}}^{\bfH}$ is a shortest representative of some coset from  $\bfW/\bfG$;
\item\label{t11.6} as objects in $\cX$, the tilting $\ttB^{\bfH}_{\bfG}$-modules are
$T_{\{e\}}^{\bfH}(w)$, where $w\in \bfW_{\{e\}}^{\bfH}$ is a shortest representative 
of some coset from $\bfW/\bfG$.
\end{enumerate}
\end{theorem}

\begin{proof}
\eqref{t11.3} follows from \cite[Section~5]{Au}. 
The functor $\theta_{\bfG}^{out}$ is exact and sends indecomposable projectives to
indecomposable projectives (see the proof of Lemma~\ref{s2.2.l201}). From the definition
of $f_{\bfG}^{\bfH}$ it follows that $\theta_{\bfG}^{out}$ sends indecomposable
$\ttA^{\bfH}_{\bfG}$-projectives to objects of the category
$\mathrm{add}\,\left(\ttA_{\{e\}}^{\bfH} f_{\bfG}^{\bfH}\right)$. This implies that 
$\theta_{\bfH}^{out}$ maps $\ttA_{\bfG}^{\bfH}\mathrm{-mod}$ to $\cX$, in particular,
the images of standard $\ttA^{\bfH}_{\bfG}$-modules belong to $\cX$. 
Furthermore, every indecomposable projective $\ttB^{\bfH}_{\bfG}$-module
has a filtration, whose subquotients are the images of standard 
$\ttA^{\bfH}_{\bfG}$-modules, since $\ttA^{\bfH}_{\bfG}$ is quasi-hereditary and 
$\theta_{\bfG}^{out}$ is exact. The existence of the duality for $\ttB^{\bfH}_{\bfG}$ 
is proved in the same way as \cite[Proposition~2.6]{MaSt2}. This 
proves \eqref{t11.2} and \eqref{t11.4}. 

$\theta_{\bfG}^{out}$ sends tilting modules to tilting modules. Hence
\eqref{t11.6} follows by tracking the highest weights.

Let $w$ be as in \eqref{t11.5}. Then $\Delta_{\bfG}^{\bfH}(w)^{\star}$ is the 
costandard $\ttA^{\bfH}_{\bfG}$-module. It is easy to see that 
$\theta_{\bfG}^{out}\Delta^{\bfH}_{\bfG}(w)^{\star}$ surjects onto the dual Verma 
module $\Delta^{\bfH}_{\{e\}}(ww_{0}^{\bfG})^{\star}$, and that the kernel of 
this surjection is generated by a projective module from $\cX$. Now \eqref{t11.5} 
follows from the fact that $\theta_{\bfG}^{on}$ sends 
$\Delta^{\bfH}_{\{e\}}(ww_{0}^{\bfG})^{\star}$ back to 
$\Delta_{\bfG}^{\bfH}(w)^{\star}$.

We are left to prove \eqref{t11.1}. First we note that $\theta_{\bfG}^{out}$ induces a
monomorphism from $\ttA_{\bfG}^{\bfH}$ to $\ttB_{\bfG}^{\bfH}$. 

From \cite[Theorem~6.1]{MaSt2} we have that the center of
$\ttA$ surjects onto $\End_{\ttA}\left(\theta_{\bfG}^{out}\Delta_{\bfG}(w_0^{\bfG})\right)$
and that the later algebra is isomorphic to $\bfC(\bfG)$. Note that
$\Delta_{\bfG}(w)\subset \Delta_{\bfG}(w_0^{\bfG})$ for all
$w\in \bfW_{\bfG}$, which induces $\theta_{\bfG}^{out}\Delta_{\bfG}(w)\subset
\theta_{\bfG}^{out}\Delta_{\bfG}(w_0^{\bfG})$.

\begin{lemma}\label{t11.lemma1}
Let $w\in \bfW_{\bfG}$. Then the restriction from $\theta_{\bfG}^{out}\Delta_{\bfG}(w_0^{\bfG})$
to $\theta_{\bfG}^{out}\Delta_{\bfG}(w)$ induces an isomorphism
\begin{displaymath}
\End_{\ttA}\left(\theta_{\bfG}^{out}\Delta_{\bfG}(w_0^{\bfG})\right)\cong
\End_{\ttA}\left(\theta_{\bfG}^{out}\Delta_{\bfG}(w)\right).
\end{displaymath}
\end{lemma}

\begin{proof}
We start with the case $w=w_0$. We have $L_{\bfG}(w_0)\cong\Delta_{\bfG}(w_0)\subset
\Delta_{\bfG}(w_0^{\bfG})$ and $[\Delta_{\bfG}(w_0^{\bfG}):L_{\bfG}(w_0)]=1$ by
\cite[Section~7]{Di}. All top subquotients of the module 
$\theta_{\bfG}^{out}\Delta_{\bfG}(w_0)$ are isomorphic to $L(w_0)$.
Since $L(w_0)$ is not a subquotient of any $\theta_{\bfG}^{out}L_{\bfG}(w)$ 
other than $\theta_{\bfG}^{out}L_{\bfG}(w_0)$, it follows that
$\theta_{\bfG}^{out}\Delta_{\bfG}(w_0)$ is stable under all endomorphisms of
$\theta_{\bfG}^{out}\Delta_{\bfG}(w_0^{\bfG})$. In particular, the restriction map
for endomorphism rings is well-defined.

Using the same arguments  as in Lemma~\ref{s2.2.l201} one shows that
$\theta_{\bfG}^{out}\Delta_{\bfG}(w_0)$ has simple top. The standard properties of 
the translation functors imply $[\theta_{\bfG}^{out}\Delta_{\bfG}(w_0):L(w_0)]=|\bfG|$,
which, together with \eqref{t11.2}, implies the equality
$\dim \End_{\ttA_{\bfG}}\left(\theta_{\bfG}^{out}\Delta_{\bfG}(w_0)\right)=|\bfG|$.
Since $\theta_{\bfG}^{out}\Delta_{\bfG}(w_0^{\bfG})$ has a Verma flag, its socle 
consists of simple subquotients isomorphic to $L(w_0)$. The arguments of the previous
paragraph imply $[\theta_{\bfG}^{out}\Delta_{\bfG}(w_0^{\bfG})/
\theta_{\bfG}^{out}\Delta_{\bfG}(w_0):L(w_0)]=0$, which, in turn, implies that the
restriction map for the endomorphism rings in injective. Since both these endomorphism rings
have the same dimension, we derive that the restriction is in fact an isomorphism. This
completes the proof in the case $w=w_0$.

In the general case it would suffice to show that the restriction from
$\theta_{\bfG}^{out}\Delta_{\bfG}(w)$ to $\theta_{\bfG}^{out}\Delta_{\bfG}(w_0)$
induces an isomorphism between the endomorphism rings of these modules. 
However, the same arguments as above show that 
$\dim \End_{\ttA_{\bfG}}\left(\theta_{\bfG}^{out}\Delta_{\bfG}(w)\right)=|\bfG|$
and that the restriction map is injective. The statement follows.
\end{proof}

\begin{lemma}\label{t11.lemma2}
Let $w\in \bfW_{\bfG}$. Then the canonical surjection
$\theta_{\bfG}^{out}\Delta_{\bfG}(w)\tto \theta_{\bfG}^{out}\Delta_{\bfG}^{\bfH}(w)$ 
induces an isomorphism,
\begin{displaymath}
\End_{\ttA}\left(\theta_{\bfG}^{out}\Delta_{\bfG}(w)\right)\cong
\End_{\ttA}\left(\theta_{\bfG}^{out}\Delta_{\bfG}^{\bfH}(w)\right).
\end{displaymath}
\end{lemma}

\begin{proof}
The arguments, analogous to those  in the proof of Lemma~\ref{t11.lemma1}, show that 
both endomorphism rings have dimension $|\bfG|$ and that the induced map is injective.
The statement follows.
\end{proof}

From Lemma~\ref{t11.lemma1} and Lemma~\ref{t11.lemma2} it follows that
\begin{displaymath}
\End_{\ccX}\left(\theta_{\bfG}^{out}\Delta_{\bfG}^{\bfH}(w)\right)=
\End_{\ttA}\left(\theta_{\bfG}^{out}\Delta_{\bfG}(w_0^{\bfG})\right)\cong
\bfC(\bfG),
\end{displaymath}
for every $w\in \bfW_{\bfG}$, moreover, that the center of 
$\ttA$ (which is $\bfC$ by \cite{So1}) surjects onto
$\End_{\ccX}\left(\theta_{\bfG}^{out}\Delta_{\bfG}^{\bfH}(w)\right)$.

Since $\bfC$ is symmetric, the corresponding 
trace form defines a splitting $\bfC(\bfG)\hookrightarrow \bfC$ of the 
epimorphism $\bfC\tto \bfC(\bfG)$ constructed above. 
This allows us to consider $\bfC(\bfG)$ as a subalgebra of $\ttB_{\bfG}^{\bfH}$, which 
is central and surjects onto all the endomorphism rings from the previous paragraph. 
The fact that $\ttB_{\bfG}^{\bfH}$ is properly stratified implies that all standard 
$\ttB_{\bfG}^{\bfH}$-modules are in fact free over
$\End_{\ccX}\left(\theta_{\bfG}^{out}\Delta_{\bfG}^{\bfH}(w_0^{\bfG})\right)\cong
\bfC(\bfG)$ (with the free basis given by any basis in the corresponding proper 
standard module). This, in particular, implies that $\ttB_{\bfG}^{\bfH}$ is a free 
$\bfC(\bfG)$-module of rank $\dim \ttA_{\bfG}^{\bfH}$. 

To complete the proof it is now enough to show that $\ttA_{\bfG}^{\bfH}$ and 
$\bfC(\bfG)$ generate $\ttB_{\bfG}^{\bfH}$. We do this using the
induction with respect to the partial order on $\bfW_{\bfG}^{\bfH}$, which makes
$\ttB_{\bfG}^{\bfH}$ stratified, see \eqref{t11.2}. Let $x,y\in \bfW_{\bfG}^{\bfH}$.
By induction we can consider a properly stratified quotient, $\tilde{\ttB}$, of
$\ttB_{\bfG}^{\bfH}$, for which $x$ becomes a maximal element. Let also $\tilde{\ttA}$
denote the corresponding quotient of $\ttA_{\bfG}^{\bfH}$. In this case we have that 
the corresponding projective $\tilde{\ttB}$-module $\tilde{P}(x)$ is standard. Then 
the properly stratified structure guarantees that the trace of $\tilde{P}(x)$ in 
the projective module $\tilde{P}(y)$ is a direct sum of say $k$ copies of 
$\tilde{P}(x)$. This implies that $\Hom_{\tilde{\ttB}}\left(\tilde{P}(x),
\tilde{P}(y)\right)$ is a free $\bfC(\bfG)$-module of rank $k$.

Let $\tilde{P}^{\tilde{\ttA}}(x)$ and $\tilde{P}^{\tilde{\ttA}}(y)$ denote the 
indecomposable projective $\tilde{\ttA}$-modules, which correspond to $x$ and $y$ 
respectively. Recall that $\theta_{\bfG}^{out}$ sends indecomposable projective
$\ttA_{\bfG}^{\bfH}$-modules to indecomposable projective $\ttB_{\bfG}^{\bfH}$-modules 
(realized as objects of $\cX$) preserving the standard filtration, and induces a 
bijection on the isomorphism classes of simple modules, compatible with the partial 
orders on these sets involved in the quasi-hereditary and properly stratified 
structures respectively. This implies that 
\begin{displaymath}
\dim\Hom_{\tilde{\ttA}}\left(\tilde{P}^{\tilde{\ttA}}(x),
\tilde{P}^{\tilde{\ttA}}(y)\right)=k,
\end{displaymath}
which shows that $\Hom_{\tilde{\ttB}}\left(\tilde{P}(x),\tilde{P}(y)\right)$ is
generated by $\ttA_{\bfG}^{\bfH}$ and $\bfC(\bfG)$. For reversed $x$ and $y$ the 
necessary statement follows applying $\star$. This proves the ungraded version of
\eqref{t11.1} and the graded version follows from the remark that all the above
arguments are compatible with the grading. This completes the proof.
\end{proof}

\begin{remark}
It follows immediately from the proof of Theorem~\ref{t11} that both standard and proper
standard $\ttB_{\bfG}^{\bfH}$-modules are gradeable, in particular, it follows that 
$\ttB_{\bfG}^{\bfH}$ is graded properly stratified in the sense of \cite[Section~8]{MaSt}.
\end{remark}


\section{Beilinson-Ginzburg-Soergel's Theorem and Backelin's Theorem}\label{s7}

\begin{theorem}\label{tbgs}(\cite[Theorem~1.1.3]{BeGiSo})
Let $\bfG$ be a parabolic subgroup of $\bfW$. Then
the algebras $\ttA_{\bfG}$ and $\ttA_{\{e\}}^{\bfG}$
are Koszul dual to each other.
\end{theorem}

\begin{proof}
Since the algebra $\bfC(\bfG)$ lives in even degrees only, from
Theorem~\ref{t11}\eqref{t11.1} it follows that
\begin{equation}\label{t22.e1}
\left(\ttB_{\bfG}^{\{e\}}\right)^!\cong \ttA_{\bfG}^!.
\end{equation}

Note again that $\ttA^!\cong\ttA$ by \cite[Theorem~18]{So1} and \cite[Theorem~2.10.1]{BeGiSo}.
Applying to the latter formula  Proposition~\ref{pdual} and using the Ringel
duality gives $\left(\ttB_{\bfG}^{\{e\}}\right)^!\cong \ttA_{\{e\}}^{\bfG}$, which, 
using \eqref{t22.e1}, implies $\ttA_{\bfG}^!\cong \ttA_{\{e\}}^{\bfG}$. But from 
Theorem~\ref{s2.2.t2} and Theorem~\ref{s5.t1} we know that both $\ttA_{\bfG}$ and 
$\ttA_{\{e\}}^{\bfG}$ are Koszul with respect to the natural grading. The proof is 
now completed by applying 
\cite[Theorem~2.10.1]{BeGiSo}.
\end{proof}

\begin{theorem}\label{tback}(\cite[Theorem~1.1]{Ba})
Let $\bfG$ and $\bfH$ be parabolic subgroups of $\bfW$. Then
the algebras $\ttA_{\bfG}^{\bfH}$ and $\ttA_{w_0 \bfH w_0}^{\bfG}$
are Koszul dual to each other.
\end{theorem}

\begin{proof}
Since the algebra $\bfC(w_0 \bfH w_0)$ lives in even degrees only, from
Theorem~\ref{t11}\eqref{t11.1} it follows that
\begin{equation}\label{t22.e101}
\left(\ttB_{w_0 \bfH w_0}^{\bfG}\right)^!\cong 
\left(\ttA_{w_0 \bfH w_0}^{\bfG}\right)^!.
\end{equation}

By Theorem~\ref{tbgs} and \cite[Theorem~2.10.1]{BeGiSo} we have
$\left(\ttA_{\{e\}}^{\bfG}\right)^!\cong \ttA^{\{e\}}_{\bfG}$.
Applying to the latter formula Proposition~\ref{pdual} and using the Ringel
duality gives $\left(\ttB_{w_0 \bfH w_0}^{\bfG}\right)^!\cong \ttA_{\bfG}^{\bfH}$, 
which, using \eqref{t22.e101}, implies $\left(\ttA_{w_0 \bfH w_0}^{\bfG}\right)^!\cong 
\ttA_{\bfG}^{\bfH}$. But from Theorem~\ref{s5.t1} we know that both 
$\ttA_{\bfG}^{\bfH}$ and $\ttA_{w_0 \bfH w_0}^{\bfG}$ are Koszul
with respect to the natural grading. The proof is now completed by applying 
\cite[Theorem~2.10.1]{BeGiSo}.
\end{proof}

\section{The category of linear complexes of tilting modules for the category 
$\ccO(\mathfrak{p},\Lambda)$}\label{s8}

Let $\bfG$ and $\bfH$ be two parabolic subgroups of $\bfW$. Let
$\bfV_{\bfG}^{\bfH}$ denote the set of all $w\in \bfW_{\bfG}$ which are at 
the same time longest coset representatives for cosets from 
$\bfH\,\backslash\bfW$. Let $q_{\bfG}^{\bfH}$ be the sum of all primitive 
idempotents in $\ttA_{\bfG}$, which correspond to $w\in \bfV_{\bfG}^{\bfH}$.
Set $\ttC_{\bfG}^{\bfH}=q_{\bfG}^{\bfH}\ttA_{\bfG} q_{\bfG}^{\bfH}$. The algebra
$\ttC_{\bfG}^{\bfH}$ is the basic associative algebra of the 
$\bfG$-singular block of the category $\OO$, studied in \cite{FuKoMa,MaSt2},
where $\mathfrak{p}$ is the parabolic subalgebra of $\fg$, associated with
$\bfH$. The algebra $\ttC_{\bfG}^{\bfH}$ is properly stratified and has
a duality. Abusing notation we denote in this section
the standard $\ttC_{\bfG}^{\bfH}$-modules by $\Delta(w)$,
$w\in \bfV_{\bfG}^{\bfH}$. The proper standard modules will be denoted by
$\overline{\Delta}(w)$, $w\in \bfV_{\bfG}^{\bfH}$.

Let now $g^{\bfG}_{\bfH}$ be the sum of all primitive idempotents in
$\ttA_{\{e\}}^{\bfG}$, which correspond to all
$w\in \bfW$ such that $w$ is the shortest coset representative for
a coset from $\bfG\backslash\bfW$ and $w$ is the shortest coset representative
for a coset from $\bfW/\bfH$ at the same time. Set 
$\ttD^{\bfG}_{\bfH}=\ttA_{\{e\}}^{\bfG}/\left(\ttA_{\{e\}}^{\bfG}
(1-g^{\bfG}_{\bfH})\ttA_{\{e\}}^{\bfG}\right)$.

$\ttC_{\bfG}^{\bfH}$ inherits a $\Z$-grading from $\ttA_{\bfG}$ and
$\ttD^{\bfG}_{\bfH}$ inherits a $\Z$-grading from $\ttA_{\{e\}}^{\bfG}$. 
We will call these gradings {\em natural}.

\begin{proposition}\label{s5.p1}
There is an isomorphism of algebras,
$\left(\ttC_{\bfG}^{\bfH}\right)^!\cong \ttD_{w_0\bfH w_0}^{\bfG}$,
compatible with the natural gradings.
\end{proposition}

\begin{proof}
This follows immediately from Proposition~\ref{pdual} and Theorem~\ref{tbgs}.
\end{proof}

In particular, we obtain that $\left(\ttC_{\{e\}}^{\bfH}\right)^!\cong 
\ttD_{w_0\bfH w_0}^{\{e\}}$ is an analogue of the algebra 
$\ttA_{\{e\}}^{w_0\bfH w_0}$ but with respect to the representatives in 
$\bfW$ for cosets on the different side.
The algebra $\ttC_{\bfG}^{\bfH}$ has both tilting and cotilting modules, 
associated with the properly stratified structure. Via \cite[Section~8]{MaSt2}
the natural grading on projective $\ttC_{\bfG}^{\bfH}$-modules induces
natural graded lifts of simple, injective, standard, proper standard, costandard 
and proper costandard modules. Following \cite[Section~5]{MaOv},
this allows us to fix a graded lift of the tilting module $T(x)$ such that
the natural inclusion $\Delta(x)\hookrightarrow T(x)$ is a homogeneous map of
degree $0$. Further, we fix the graded lift of the cotilting module $C(x)$ 
such that the natural projection $C(x)\tto\nabla(x)$ is a homogeneous map of
degree $0$. Note that $T(x)\cong C(x)\langle -2\mathfrak{l}(w_0^{\bfH})\rangle$.

\begin{theorem}\label{s5.t2}
\begin{enumerate}[(i)]
\item\label{s5.t2.1} All standard $\ttC_{\{e\}}^{\bfH}$-modules admit linear 
tilting coresolutions.
\item\label{s5.t2.2} All costandard $\ttC_{\{e\}}^{\bfH}$-modules admit linear 
cotilting resolutions.
\end{enumerate}
\end{theorem}

\begin{proof}
Because of the duality it is enough to prove \eqref{s5.t2.1}.
In this case we have $\Delta(w_0)=T(w_0)$ for which the statement
is trivial. Now the statement follows by applying induction on the length of $w$ 
and using translation functors, and the same arguments as in the proof of 
Proposition~\ref{s2p1}. 
\end{proof}

\begin{remark}
One can also prove Theorem~\ref{s5.t2} generalizing the arguments of 
Proposition~\ref{s2.2.l2}. The main difference is that one will be forced to
write an analogue of the parabolically induced BGG-resolution for cosets 
$\bfG/\bfH$. However, since this resolution consists only of modules with
the scalar action of the center of $\mathfrak{g}$, the existence of such 
resolution follows directly from the proof of Proposition~\ref{s2.2.l2} using 
the functor $\eta$ from \cite[Section~3]{MaSt2}.
\end{remark}

We believe that the statement of Theorem~\ref{s5.t2} remains valid for singular
blocks (that is for all $\ttC_{\bfG}^{\bfH}$) as well. However, the translation 
functor approach does not work appropriately for singular blocks, and a proper 
generalization of the corresponding arguments of Theorem~\ref{s2.2.t2} seems to 
be technically very complicated. 

\begin{corollary}\label{s5.c3}
\begin{enumerate}[(i)]
\item\label{s5.c3.1} All standard $\ttC_{\{e\}}^{\bfH}$-modules admit linear 
projective resolutions.
\item\label{s5.c3.2} All costandard $\ttC_{\{e\}}^{\bfH}$-modules admit linear 
injective coresolutions.
\end{enumerate}
\end{corollary}

\begin{proof}
This follows from Theorem~\ref{s5.t2} and the Ringel self-duality of 
$\ttC_{\{e\}}^{\bfH}$, see \cite[Theorem~3]{FuKoMa}.
\end{proof}

\begin{corollary}\label{s5.c4}
Let $x,y\in \bfV_{\{e\}}^{\bfH}$ and $x>y$. Then
\begin{displaymath}
\ext_{\ttC_{\{e\}}^{\bfH}}^k
\left(\overline{\Delta}(x)\langle m\rangle,\Delta(y)\right)\neq 0
\end{displaymath}
implies $m+2\mathfrak{l}(w_0^{\bfH})\leq k\leq \mathfrak{l}(x)-\mathfrak{l}(y)$. 
In particular,
\begin{displaymath}
\Ext_{\ttC_{\{e\}}^{\bfH}}^k\left(\overline{\Delta}(x),\Delta(y)\right)=
\Ext_{\ttC_{\{e\}}^{\bfH}}^k\left(\Delta(x),\Delta(y)\right)=0
\end{displaymath}
for all $k\geq \mathfrak{l}(x)-\mathfrak{l}(y)$.
\end{corollary}

\begin{proof}
Since all tilting $\ttC_{\{e\}}^{\bfH}$-modules are cotilting at the same time,
see for example \cite[Section~6]{FuKoMa}, 
the tilting coresolution $\mathcal{T}^{\bullet}(\Delta(y))$ of $\Delta(y)$ 
is cotilting (up to a shift of grading) at the same time. But this means that 
\begin{displaymath}
\ext_{\ttC_{\{e\}}^{\bfH}}^k\left(\overline{\Delta}(x)\langle m\rangle,\Delta(y)\right)=
\mathrm{Hom}(\overline{\Delta}(x)\langle l\rangle^{\bullet},
\mathcal{T}^{\bullet}(\Delta(y))),
\end{displaymath}
where the last homspace is taken in the homotopy category of the category
$\ttC_{\{e\}}^{\bfH}\mathrm{-grmod}$. Now the first statement follows from
Theorem~\ref{s5.t2}\eqref{s5.t2.1}, Proposition~\ref{p5}, and the above remark
that $T(x)\cong C(x)\langle -2\mathfrak{l}(w_0^{\bfG})\rangle$. The second 
statement follows from the first one.
\end{proof}

\begin{corollary}\label{s5.c5}
Let $x,y\in \bfV_{\{e\}}^{\bfH}$ and $x>y$. Then
\begin{displaymath}
\Ext_{\ttC_{\{e\}}^{\bfH}}^{\mathfrak{l}(x)-\mathfrak{l}(y)}
\left(\Delta(x),\Delta(y)\right)
\end{displaymath}
is a free $\mathrm{End}_{\ttC_{\{e\}}^{\bfH}}(\Delta(x))$-module of rank
\begin{displaymath}
\dim\ext_{\ttC_{\{e\}}^{\bfH}}^{\mathfrak{l}(x)-\mathfrak{l}(y)}
\left(\overline{\Delta}(x)
\langle \mathfrak{l}(x)-\mathfrak{l}(y)-2\mathfrak{l}(w_0^{\bfH})\rangle,\Delta(y)\right).
\end{displaymath}
\end{corollary}

\begin{proof}
As in Corollary~\ref{s5.c4} we have that 
\begin{displaymath}
\Ext_{\ttC_{\{e\}}^{\bfH}}^{\mathfrak{l}(x)-\mathfrak{l}(y)}\left(\Delta(x),\Delta(y)\right)=
\mathrm{Hom}(\Delta(x)^{\bullet},\mathcal{T}^{\bullet}(\Delta(y))),
\end{displaymath}
where the last homspace is taken in the homotopy category of the category
$\ttC_{\{e\}}^{\bfH}\mathrm{-mod}$. Note that all (graded) standard 
$\ttC_{\{e\}}^{\bfH}$-modules have a graded filtration by proper standard 
$\ttC_{\{e\}}^{\bfH}$-modules, see \cite[Section~8]{MaSt2}, and further observe that 
\begin{displaymath}
\dim\ext_{\ttC_{\{e\}}^{\bfH}}^{\mathfrak{l}(x)-\mathfrak{l}(y)}\left(\overline{\Delta}(x)
\langle \mathfrak{l}(x)-\mathfrak{l}(y)-2\mathfrak{l}(w_0^{\bfH})\rangle,\Delta(y)\right)
\end{displaymath}
is the number of occurrences of $T(x)\langle \mathfrak{l}(x)-\mathfrak{l}(y)\rangle$
in $\mathcal{T}^{\mathfrak{l}(x)-\mathfrak{l}(y)}(\Delta(y))$. Further, 
$\mathrm{Hom}_{\ttC_{\{e\}}^{\bfH}}(\Delta(x),T(x))\cong 
\mathrm{Hom}_{\ttC_{\{e\}}^{\bfH}}(\Delta(x),\Delta(x))$, which is a free
$\mathrm{End}_{\ttC_{\{e\}}^{\bfH}}(\Delta(x))$-module of rank one. The statement follows.
\end{proof}

Though $\ttD_{\bfG}^{\bfH}$ is not quasi-hereditary in general (a counter example can be
derived from \cite[Remark~1.2]{MaSt}), the notion of a standard module is nevertheless
well-defined for this algebra (and for any algebra with a fixed order on the set of 
the isomorphism classes of simple modules). 

\begin{corollary}\label{s5.c7}
Under the equivalence provided by Proposition~\ref{s5.p1},
the linear tilting coresolutions of standard $\ttC_{\{e\}}^{\bfH}$-modules are 
standard $\ttD_{w_0\bfH w_0}^{\{e\}}$-modules; and (appropriately shifted) 
linear cotilting resolutions of costandard $\ttC_{\{e\}}^{\bfH}$-modules are 
costandard $\ttD_{w_0\bfH w_0}^{\{e\}}$-modules.
\end{corollary}

\begin{proof}
The proof is analogous to that of \cite[Proposition~5]{MaOv}.
\end{proof}

\begin{remark}
All standard $\ttC_{\{e\}}^{\bfH}$-modules are quotients of Verma modules in $\cO$ 
and hence have central character. Applying the equivalence $\eta$ from
\cite[Section~3]{MaSt2}, we obtain  that the multiplicities of simple modules 
in the composition series of standard $\ttC_{\{e\}}^{\bfH}$-modules
coincide with the multiplicities of the corresponding (under $\eta$) simple 
$\ttA_{\{e\}}^{\bfH}$-modules in the composition series of the corresponding (again
under $\eta$) standard $\ttA_{\{e\}}^{\bfH}$-modules. Moreover, everything is 
compatible with the grading. Therefore, the components of the linear tilting 
coresolution of a standard $\ttC_{\{e\}}^{\bfH}$-modules can be computed using the 
Kazhdan-Lusztig combinatorics (see also \cite{CC}).
\end{remark}

\vspace{0.2cm}

\begin{center}
\bf Acknowledgments
\end{center}

The research was partially supported by the Swedish Research Council, the
Royal Swedish Academy of Sciences, the Swedish Foundation for International Cooperation
in Research and Higher Education (STINT), and by the Faculty of Natural Sciences,
Uppsala University. I would like to thank Yuriy Drozd and Ryszard Rubinsztein for 
helpful and stimulating discussions. I am especially grateful to Catharina Stroppel
for many helpful discussions, suggestions and corrections to the preliminary version of 
the paper.

\vspace{0.5cm}

\noindent
Volodymyr Mazorchuk, Department of Mathematics, Uppsala University,
Box 480, 751 06, Uppsala, Sweden,
e-mail: {\tt mazor\symbol{64}math.uu.se},\\
web: {``http://www.math.uu.se/$\tilde{\hspace{1mm}}$mazor/''}.

\end{document}